\documentclass[11pt]{article}

\usepackage{amsmath}
\usepackage{amssymb}
\usepackage{graphicx}
\usepackage{url}
\usepackage{a4}

\def\cqfd{\skip10=\parfillskip\parfillskip=0pt
\enspace\hfill\symbolecqfd\par\parfillskip=\skip10\par\medskip}
\def\symbolecqfd{\rlap{$\sqcap$}$\sqcup$}

\newtheorem{theorem}{Theorem}[section]
\newtheorem{proposition}[theorem]{Proposition}
\newtheorem{lemma}[theorem]{Lemma}
\newtheorem{corollary}[theorem]{Corollary}

\newtheorem{pro-fact}[theorem]{Fact}
\newenvironment{fact}{\begin{pro-fact}\rm}{\cqfd\end{pro-fact}}

\newtheorem{pro-example}[theorem]{Example}

\newtheorem{pro-remark}[theorem]{Remark}
\newenvironment{remark}{\begin{pro-remark}\rm}{\cqfd\end{pro-remark}}

\newenvironment{preuve}{\rm \trivlist \item[\hskip \labelsep{\bf
Proof.}]}{\cqfd\endtrivlist}

\def\cqfd{\skip10=\parfillskip\parfillskip=0pt
\enspace\hfill\symbolecqfd\par\parfillskip=\skip10\par\medskip}
\def\symbolecqfd{\rlap{$\sqcap$}$\sqcup$}
\def\proof{\begin{preuve}}
\def\eop{\end{preuve}}

\def\proofof#1{\rm \trivlist \item[\hskip \labelsep{\bf
Proof of~#1.}]}
\def\eopo{\cqfd\endtrivlist}

\newcounter{commentcounter}

\def\rank{\textsf{rank}}
\def\rrank{\widetilde{\textsf{rk}}}

\def\inter[#1]{[\![#1]\!]}
\def\inv{^{-1}}

\def \O {\mathcal{O}}

\def \calS {\mathcal{S}}
\def \calP {\mathcal{P}}
\def \calF {\mathcal{F}}
\def \calJ {\mathcal{J}}
\def \calK {\mathcal{K}}
\def \calL {\mathcal{L}}

\let\phi\varphi
\let\epsilon\varepsilon

\def\Pref{\textsf{Pref}}

\begin{document}

\title{Statistical properties of subgroups of free groups\thanks{%
All authors benefitted from the support of the French-Spanish program
\textsc{picasso} (project ACI-HF2006-0239).  The first and third
authors were supported by the ANR (projects \textsc{blan 07-2\_195422 gamma} and \textsc{anr-2010-blan-0204 magnum}).
The second and fourth authors were supported by MEC (Spain) and the
EFRD (EC) through project number MTM2008-01550.  The last author was
supported by the ESF program \textsc{AutoMathA} and by the ANR (project \textsc{anr-2010-blan-0202 frec}).}%
}

\markboth{}{}

\author{
    Fr\'ed\'erique Bassino, \small{\url{Frederique.Bassino@lipn.univ-paris13.fr}}\\
    \small{Universit\'e Paris 13, LIPN}\thanks{%
    LIPN, Universit\'e Paris 13 and CNRS, Institut Galil\'ee, 99,
    avenue Jean-Baptiste Cl\'ement, 93430 Villetaneuse, France}
    \and
    Armando Martino, \small{\url{A.Martino@upc.soton.ac.uk}}\\
    \small{University of Southampton}\thanks{%
    School of Mathematics, University of Southampton, University Road,
    Southampton SO17 1BJ, United Kingdom}
    \and
    Cyril Nicaud, \small{\url{nicaud@univ-mlv.fr}}\\
    \small{Universit\'e Paris-Est, LIGM and CNRS}\thanks{%
    Institut Gaspard Monge, Universit\'e Paris-Est, 77454
    Marne-la-Vall\'ee Cedex 2, France}
    \and
    Enric Ventura, \small{\url{enric.ventura@upc.edu}}\\
    \small{Universitat Politecnica de Catalunya}\thanks{%
    Universitat Politecnica de Catalunya, Av.  Bases de Manresa 61-73,
    08240 Manresa, Catalunya, Spain }
    \and
    Pascal Weil, \small{\url{pascal.weil@labri.fr}}\\
    \small{Univ. Bordeaux, LaBRI, UMR 5800, F-33400 Talence, France}\\
    \small{CNRS, LaBRI, UMR 5800, F-33400 Talence, France}\thanks{%
    LaBRI, Univ. Bordeaux, 351 cours de la Lib\'eration, 33400 Talence, France.}
    }

    \date{}

\maketitle

\newpage

\begin{abstract}
    The usual way to investigate the statistical properties of
    finitely generated subgroups of free groups, and of finite
    presentations of groups, is based on the so-called word-based
    distribution: subgroups are generated (finite presentations are
    determined) by randomly chosen $k$-tuples of reduced words, whose
    maximal length is allowed to tend to infinity.  In this paper we
    adopt a different, though equally natural point of view: we
    investigate the statistical properties of the same objects, but
    with respect to the so-called graph-based distribution, recently
    introduced by Bassino, Nicaud and Weil.  Here, subgroups (and
    finite presentations) are determined by randomly chosen Stallings
    graphs whose number of vertices tends to infinity.

    Our results show that these two distributions behave quite
    differently from each other, shedding a new light on which
    properties of finitely generated subgroups can be considered
    \textit{frequent} or \textit{rare}.  For example, we show that
    malnormal subgroups of a free group are negligible in the
    graph-based distribution, while they are exponentially generic in
    the word-based distribution.  Quite surprisingly, a random finite
    presentation generically presents the trivial group in this new
    distribution, while in the classical one it is known to
    generically present an infinite hyperbolic group.
\end{abstract}

\bigskip


\medskip


\bigskip


\section{Introduction}

Statistical properties of elements and subgroups of free groups have
evoked much interest in recent years, especially after Gromov's famous
claim \cite[0.2.A]{1987Gromov} that ``most" groups were hyperbolic, which led to precise statements
and proofs by Ol'shanski\u\i\ \cite{Olshanskii1992} and by Champetier \cite{1991Champetier,1995Champetier}.
Shortly thereafter, Ol'shanski\u\i\ and Arzhantseva
\cite{1996ArjantsevaO,Arjantseva2000} pursued the study of the
statistical properties of finite presentations of groups, that is,
largely, of finitely generated normal subgroups of free groups. We refer the reader to the survey by Ollivier \cite{2005Ollivier} for more details.

This interest encountered another historical trend in combinatorial
group theory, namely the consideration of algorithmic problems, which
leads naturally to an interest in the evaluation of the complexity of
these algorithms (e.g. \cite{2000BirgetMMW,2001MargolisSW,2007RoigVW})
and in enumeration problems.

The search for innovative group-based cryptographic systems (see
\cite{2008bookCRM} for instance) only reinforced the study of complexity
questions, and focused it on the investigation of the statistical
properties of finitely generated subgroups of free groups, notably via
the notion of generic complexity (see~\cite{2006KSS,2002Jitsukawa}).

The usual method to approach statistical properties is to enumerate
the objects under consideration, or more precisely, {\em
representatives} for these objects, in a stratified way.  For
instance, if we wish to investigate $k$-generated subgroups of $F_r$
(resp.  finitely presented groups with $r$ generators and $k$
relators), we proceed by enumerating lists of $k$-tuples of generators
(resp.  relators) over a fixed alphabet of $r$ letters, so that at
level $n$ one has enumerated all such $k$-tuples whose elements have
length at most $n$.  In the situation we will consider, there are only
finitely many objects of a given level $n$ and it makes sense to ask
what proportion of level $n$ objects satisfy a given property.  This
gives us a number $p_n$ between 0 and 1 for each $n$, associated to
the given property and one can ask whether this sequence has a well
defined limit.  If the limit exists and is equal to $1$, we would say
that the property is {\em generic}, and take this to mean that most
objects satisfy the property.  At the other extreme, if the limit of
the $p_n$ equals 0, we would say that the property in question is {\em
negligible} and conclude that it is rarely encountered amongst our
objects.

A crucial observation, which is well worth mentioning in view of the
intuitive weight carried by expressions such as \textit{most objects}
or \textit{rarely encountered}, is that genericity and negligibility
depend essentially on the choice of the stratification: different
stratifications of the same objects, say finitely generated subgroups
of free groups, will bring to light different insights on the
statistical behavior of these objects.  Concretely, different
properties will appear to be generic or negligible.

Up to recently (namely the publication of \cite{BNW}), the literature
was unanimous in adopting the representation of finitely generated
subgroups of free groups by $k$-tuples of generators, stratified by
their maximal length -- which we call the word-based distribution.

It is the purpose of this paper to question this unanimity. The basic 
idea is that there exists another very natural representation of 
finitely generated subgroups of free groups, by their Stallings graph 
(\cite{1983Stallings}, see Section~\ref{sec: stallings}). Stratifying 
finitely generated subgroups by the size (number of vertices) of 
their Stallings graph -- what we call the graph-based distribution -- 
indeed sheds a different light on which properties of subgroups are 
frequent or rare. One of our main results is that malnormality and 
purity, which are generic in the word-based distribution, are 
negligible in the graph-based distribution (Section~\ref{sec: 
malnormal}).

We also exhibit a property of finitely generated subgroups of $F_r$ 
that is negligible in the word-based distribution and that has a 
non-zero, non-one asymptotic probability (namely $e^{-r}$) in the 
graph-based distribution (Section~\ref{sec: intermediate}).

Finally we explore the possibility of using the graph-based
distribution to discuss the statistical properties of finitely
presented groups.  The results there are disappointing: it turns out
that finitely presented groups are generically trivial in this
distribution -- quite differently from the word-based distribution in
which they are known to be generically infinite and hyperbolic
(Section~\ref{sec: fp groups}).

Sections~\ref{sec: prelim} and~\ref{sec: 2 distributions} are devoted 
to preliminaries on genericity and to a review of the main features 
of the word-based and the graph-based distributions for finitely 
generated subgroups of free groups.

\section{Preliminaries}\label{sec: prelim}

Here we summarize standard facts about the Stallings graphs of
subgroups (in Section~\ref{sec: stallings}) and we review the notions
of generic and negligible properties.

Throughout the paper, $A$ denotes an alphabet, that is, a finite
non-empty set and $F(A)$ denotes the free group over $A$.  The
elements of $F(A)$ are represented by the reduced words written using
letters from $A$ and their formal inverses $\{a\inv \mid a\in A\}$.
If $r\ge 1$, we often use the notation $F_r$ instead of $F(A)$, to
indicate that $A$ consists of $r$ letters.  Throughout the paper, we
will in fact assume that $r\ge 2$.

We denote by $[n]$ ($n\ge 1$) the set $\{1,\ldots, n\}$.

\subsection{Subgroup graph representation}\label{sec: stallings}

Each finitely generated subgroup of $F(A)$ can be represented uniquely
by a finite graph of a particular type, by means of the technique
known as \textit{Stallings foldings} \cite{1983Stallings}. This representation has been used by many authors, frequently using combinatorial, graph-theoretic notations that slightly differ from those used by Stallings. It is this formalism that we also use, which can also be found in
\cite{2000Weil,2002KapovichMyasnikov,2006Touikan,2007MiasnikovVW}.
The procedure of Stallings foldings is informally described at the end of this section.
 
An \textit{$A$-graph} is defined to be a pair $\Gamma = (V,E)$ with $E
\subseteq V\times A\times V$, such that
\begin{itemize}
    \item if $(u,a,v), (u,a,v') \in E$, then $v = v'$;
    \item if $(u,a,v), (u',a,v) \in E$, then $u = u'$.
\end{itemize}
The elements of $V$ are called the \textit{vertices} of $\Gamma$ and
the elements of $E$ are its \textit{edges}.  We say that $\Gamma$ is
\textit{connected} if the underlying undirected graph is connected.
If $v\in V$, we say that $v$ is a \textit{leaf} if $v$ occurs at most
once in (the list of triples defining) $E$ and we say that $\Gamma$ is
\textit{$v$-trim} if no vertex $w\ne v$ is a leaf.  Finally we say
that the pair $(\Gamma,v)$ is \textit{admissible} if $\Gamma$ is a 
finite,
$v$-trim and connected $A$-graph.  Then it is known (see
\cite{1983Stallings,2000Weil,2002KapovichMyasnikov,2007MiasnikovVW})
that:
\begin{itemize}
  \item Stallings  associated with each finitely generated
    subgroup $H$ of $F(A)$ a unique admissible pair of the form
    $(\Gamma,1)$, which we call the \textit{graphical representation} 
    or the \textit{Stallings graph} of $H$ and write $\Gamma(H)$;
    
    \item every admissible pair $(\Gamma,1)$ is the graphical
    representation of a unique finitely generated subgroup of $F(A)$;
    
    \item if $(\Gamma,1)$ is the graphical representation of $H$ and
    $u$ is a reduced word, then $u\in H$ if and only if $u$ labels a
    loop at 1 in $\Gamma$ (by convention, an edge $(u,a,v)$ can be read from $u$ to $v$ with label $a$, or from $v$ to $u$ with label $a\inv$);
    
    \item if $(\Gamma,1)$ is the graphical representation of $H$, then
    $\rank(H) = |E| -|V| + 1$;
    
    \item finitely generated subgroups $H$ and $K$ are conjugates if
    and only if the cyclic cores of $\Gamma(H)$ and $\Gamma(K)$
    (obtained by repeatedly deleting leaves and the edges they are
    adjacent to) are equal.
\end{itemize}

We informally remind the readers of the computation of the graphical
representation of a subgroup generated by a subset $B=\{u_1,\ldots,
u_k\}$.  It consists in building an $(A\sqcup A\inv)$-graph, changing
it into a $A$-graph, then reducing it using \textit{foldings}.  First
build a vertex $1$.  Then, for every word $u$ of length $n$ in $B$,
build a loop with label $u$ from $1$ to $1$, adding $n-1$ vertices.
Change every edge $(u,a\inv,v)$ labeled by a letter of $A\inv$ into an
edge $(v,a,u)$.  Then iteratively identify the vertices $v$ and $w$
whenever there exists a vertex $u$ and a letter $a\in A$ such that
either both $(u,a,v)$ and $(u,a,w)$ or both $(v,a,u)$ and $(w,a,u)$
are edges in the graph (the corresponding two edges are
\textit{folded}, in Stallings' terminology).

The resulting graph $\Gamma$ is such that $(\Gamma,1)$ is admissible
and, very much like in the (1-dimensional) reduction of words, it does
not depend on the order used to perform the foldings.

\subsection{Negligibility and genericity}

Let $S$ be a countable set, the disjoint union of finite sets $S_n$
($n\ge 0$), and let $B_n = \bigcup_{i\le n}S_i$.  Typically in this
paper, $S$ will be the set of Stallings graphs, of partial injections,
of reduced words or of $k$-tuples of reduced words, and $S_n$ will be
the set of elements of $S$ of size $n$.
A subset $X$ of $S$ is \textit{negligible} (resp.  \textit{generic})
if the probability for an element of $B_n$ to be in $X$, tends to 0
(resp.  to 1) when $n$ tends to infinity; that is, if
$\lim_n\frac{|X\cap B_n|}{|B_n|} = 0$ (resp.  $=1$).

Naturally, the negligibility or the genericity of a subset $X$ of $S$
depends on the layering of $S$ into the $S_n$.  In particular, if $X$
and its complement are both infinite, then an appropriate partition of
$S$ into finite subsets $S_n$ will make $X$ negligible, another will
make it generic, and indeed, another will be such that
$\lim_n\frac{|X\cap B_n|}{|B_n|} = p$ for any fixed $0 < p < 1$.
    
Thus, any discussion of negligibility or genericity must clearly
specify the distribution that is considered, that is, the choice of
the partition $(S_n)_n$.

\subsubsection{Rate of convergence}
In general, we may be interested in the speed of convergence of
$\frac{|X\cap B_n|}{|B_n|}$ -- towards 0 if $X$ is negligible and
towards 1 if it is generic.  One reason is that a higher speed of
convergence indicates a higher rate of confidence that a randomly
chosen element of $S$ of size $n$ will miss $X$ if $X$ is negligible,
or will be in $X$ if $X$ is generic, even for moderately large values
of $n$.

If a class $\calF$ of functions tending to 0 is closed under $\max$
(of two elements), we say that a subset $X$ is
\textit{$\calF$-negligible} if $\frac{|X\cap B_n|}{|B_n|} = \O(f(n))$
for some $f\in \calF$.  We also say that $X$ is
\textit{$\calF$-generic} if the complement of $X$ is
$\calF$-negligible.  Note that $\calF$-negligible (resp.
$\calF$-generic) sets are closed under finite unions and
intersections.

Much of the literature is concerned with \textit{exponential}
negligibility or genericity, namely $\calF$-negligibility or
genericity where $\calF$ is the class of functions $e^{-cn}$ ($c>0$).

\subsubsection{Balls versus spheres}\label{sec: balls spheres}

The definition of negligibility and genericity above is given in terms
of the balls $B_n$: the sets of elements of size at most $n$.  It is
sometimes more expedient to reason in terms of the proportion of
elements of $X$ in the spheres $S_n$: let us say, within the ambit of
this section, that a set $X$ is S-negligible (resp.  S-generic) if the
ratio $\frac{|X\cap S_n|}{|S_n|}$ tends to $0$ (resp.  $1$).  The
definition of $\calF$-S-negligibility or $\calF$-S-genericity is
analogous.  We verify in this section that (exponential)
S-negligibility implies (exponential) negligibility.  The same holds
of course for genericity.

\begin{proposition}\label{S and B negligibility}
    An S-negligible (resp.  S-generic) set is also negligible (resp.
    generic).
    
    If the structures under consideration grow fast enough, so that
    $\lim\frac{B_n}{B_{2n}} = 0$, then the same result holds for
    exponential negligibility and genericity.
\end{proposition}

The proof of this statement relies on the following technical lemma.

\begin{lemma}\label{lemma S and B}
    Let $(a_n)$ and $(b_n)$ be increasing sequences of positive real 
    numbers.
    \begin{itemize}
	\item[(1)] (\textbf{Stolz-Ces\`aro theorem}) If $\lim b_n = \infty$ and
	$\lim\frac{a_{n+1}-a_n}{b_{n+1}-b_n} = 0$, then
	$\lim\frac{a_n}{b_n} = 0$.

	\item[(2)] If $(\frac{b_n}{b_{2n}})$ and
	$(\frac{a_{n+1}-a_n}{b_{n+1}-b_n})$ converge to $0$
	exponentially fast and if $a_n \le b_n$ for each $n$, then
	$(\frac{a_n}{b_n})$ converges to $0$ exponentially fast as
	well.
    \end{itemize}
\end{lemma}

\proof
(1)\enspace Since $\lim\frac{a_{n+1}-a_n}{b_{n+1}-b_n} = 0$, for each
$\epsilon > 0$, there exists $n_0$ such that $a_{n+1}-a_n \le
\epsilon(b_{n+1}-b_n)$ for all $n \ge n_0$.  Summing these inequalities
for all integers between $n-1$ and $n_0$, we find that $a_{n}-a_{n_0}
\le \epsilon(b_{n}-b_{n_0})$ for all $n> n_0$.  Dividing by $b_n$ and
using the fact that $\lim b_n = \infty$, we conclude that
$\frac{a_n}{b_n} < 2\epsilon$ for all large enough $n$.

(2)\enspace Our hypothesis is now that there exists $c> 0$ such that
$a_{n+1}-a_n \le e^{-cn}(b_{n+1}-b_n)$ for all $n \ge n_0$.  Summing
these inequalities for the integers between $n$ and $2n-1$, we find
that $a_{2n}-a_n \le e^{-cn}(b_{2n}-b_n)$ for all $n \ge n_0$.  We now
divide both sides by $b_{2n}$ and use the fact that
$\frac{a_n}{b_{2n}} \le \frac{b_n}{b_{2n}}$ and that this sequence
converges to $0$ exponentially fast to conclude that
$(\frac{a_{2n}}{b_{2n}})_n$ converges to $0$ exponentially fast.
Summing instead for the integers between $n$ and $2n$ and dividing by
$b_{2n+1}$ shows that $(\frac{a_{2n+1}}{b_{2n+1}})_n$ converges to $0$
exponentially fast as well.
\eop

\proofof{Proposition~\ref{S and B negligibility}}
Let $X \subseteq S$, $a_n = |X \cap B_n|$ and $b_n = |B_n|$.  Then
$a_n-a_{n-1} = |X\cap S_n|$ and $b_n - b_{n-1} = |S_n|$.  The
statement on (exponential) negligibility now follows directly from
Lemma~\ref{lemma S and B}.  The statement on genericity follows as
well, since generic sets are the complements of negligible sets.
\eopo

\section{The word-based and the graph-based distributions}
\label{sec: 2 distributions}

In order to discuss the distribution of finitely generated subgroups
of $F_r$, we need to fix a representation of these subgroups by means
of discrete structures.  In this paper we consider two such
structures: a subgroup can be given by a tuple of generators (reduced
words in $F_r$), or by its Stallings graph (Section~\ref{sec:
stallings}).  In the first case, the size of the representation is the
pair $(k,n)$ where $k$ is the number of generators and $n$ their
maximal length -- or $n$ if $k$ is fixed; in the second case, the size
of the representation is the number $n$ of vertices of the Stallings
graph.  In either case, there are only finitely many subgroups of each
size.

We first review the literature on the word-based and the graph-based
distributions (Sections~\ref{sec: word based} and~\ref{sec: graph
based}), and then start the discussion of negligible or generic
properties of subgroups (Section~\ref{sec: begin discussion}).

\subsection{The word-based distribution}\label{sec: word based}

The distribution usually found in the literature (e.g.
\cite{2006KSS,2002Jitsukawa,2003KapovichMSS}) is in fact a
distribution on the $k$-tuples $\vec h = (h_1,\ldots,h_k)$ of reduced
words of length at most $n$, where $k$ is fixed and $n$ is allowed to
grow to infinity; one then considers the subgroup $H$ generated by
$\vec h$.  We call this distribution \textit{word-based}.

Let us first record three elementary facts, which can also be found in
\cite{2002Jitsukawa}\footnote{%
We choose to reiterate the proofs of these results, because we feel
that our presentation exhibits more clearly their combinatorial
underpinnings.}%
. We denote by $R_n$ the set of reduced words of length at most $n$.

\begin{fact}\label{fact: Rn}
    ${\displaystyle|R_n| = \frac{r}{r-1}(2r-1)^n \left(1 -
    \frac1{r(2r-1)^n}\right)}$.
\end{fact}

\proof
The number of reduced words of length $i\ge 1$ is $2r(2r-1)^{i-1}$, so
the cardinality of $R_n$ is
\begin{align*}
    |R_n| = 1 + \sum_{i=1}^n2r (2r-1)^{i-1} &= 1 + 2r \frac{(2r-1)^n
    - 1}{2r-2} \\
    &= 1 + \frac{r}{r-1}((2r-1)^n - 1)\\
    &= \frac{r}{r-1}(2r-1)^n \left(1 - \frac1{r(2r-1)^n}\right).
\end{align*}
\eop

\begin{fact}\label{fact: min length}
    Let $0 < \alpha < 1$.  Exponentially generically, a reduced word
    in $R_n$ has length greater than $\alpha n$.
\end{fact}

\proof
The proportion of words in $R_n$, of length less than or equal to
$\alpha n$, is
\begin{align*}
    \frac{|R_{\lfloor \alpha n\rfloor}|}{|R_n|} =
    \frac{\frac{r}{r-1}(2r-1)^{\lfloor \alpha n\rfloor}(1 +
    o(1))}{\frac{r}{r-1}(2r-1)^n(1+o(1))} &= (2r-1)^{\lfloor \alpha
    n\rfloor - n}(1+o(1))\\
    &\le (2r-1)^{(\alpha -1)n}(1+o(1)).
\end{align*}
Since $\alpha - 1 < 0$, it converges to 0 exponentially
fast.
\eop

Let $\vec h = (h_1,\ldots,h_k)$ be a tuple of reduced words and let
$\mu>0$ such that $\min|h_i| > 2\lceil\mu\rceil$.  We denote by
$\Pref_\mu(\vec h)$ the set of prefixes of length at most
$\lceil\mu\rceil$ of the $h_i$ and $h_i\inv$.

\begin{fact}\label{fact: distinct prefixes}
    Let $0 < \lambda < \frac12$.  Exponentially generically, a
    $k$-tuple $(h_1,\ldots,h_k)$ of elements of $R_n$, is such that
    $\min|h_i| > 2\lceil\lambda n\rceil$ and the prefixes of the $h_i$
    and $h_i\inv$ of length $\lceil\lambda n\rceil$ are pairwise
    distinct.
\end{fact}

\proof
The complement in $R_n^k$ of the set of $k$-tuples described in the
statement is the union of the set $Y$ of $k$-tuples $\vec h$ where
$\min|h_i| \le 2\lceil\lambda n\rceil$, and of the set $Z$ of
$k$-tuples where $\min|h_i| > 2\lceil\lambda n\rceil$ and the set of
prefixes of length $\lceil\lambda n\rceil$ of the $h_i$ and $h_i\inv$
has at most $2k-1$ elements.  Since $2\lambda < 1$, the set $Y$ is
exponentially negligible by Fact~\ref{fact: min length} and we now
concentrate on $Z$.

For each integer $2\lceil\lambda n\rceil < m \le n$, let $Z_m$ be the
set of $k$-tuples in $Z$, such that $\min|h_i| = m$. Then
\begin{align*}
    |Z_m| &\le (2r(2r-1)^{\lceil\lambda
    n\rceil-1})^{2k-1}\ k(2k-1)\ (2r-1)^{k(m-2\lceil\lambda n\rceil)} \\
    &\le (2r)^{2k-1}\ k(2k-1)\ (2r-1)^{k(m-2)- \lceil\lambda n\rceil +1}.
\end{align*}
Summing these inequalities for all $2\lceil\lambda n\rceil < m \le 
n$, we find
$$|Z| \enspace\le\enspace (2r)^{2k-1}\ k(2k-1)\ (2r-1)^{k(n-1)-
\lceil\lambda n\rceil +1}.$$
As a result, the proportion of $k$-tuples in $Z$  is at most
$$\frac{(2r)^{2k-1}\ k(2k-1)\ (2r-1)^{k(n-1)- \lceil\lambda n\rceil
+1}}{\frac{r^k}{(r-1)^k}(2r-1)^{kn}(1+o(1))} \enspace\le\enspace
C\ (2r-1)^{-\lambda n}(1+o(1))$$
for some constant $C$ depending only on $k$ and $r$.  Thus, this
proportion converges to 0 exponentially fast.
\eop

\begin{remark}
    The small cancellation property, a closely related statement relative to common factors located
    anywhere in the words $h_i$ and $h_i\inv$ (not just at their
    extremities), is discussed in Lemma~\ref{fact: C' single word}, in
    a variant of Arzhantseva and Ol'shanski\u\i's result on cyclic
    words \cite[Lemma 3]{1996ArjantsevaO}.
\end{remark}

Let $0<\alpha<1$ and $0< \lambda<\frac\alpha 2$, and for each $n$, let
$Y_{\alpha,\lambda,n,k}$ be the set of $k$-tuples $\vec h =
(h_1,\ldots,h_k) \in R_n^k$ such that $\min|h_i| > \alpha n$ and the
prefixes of the $h_i$ and $h_i\inv$ of length $\lceil\lambda n\rceil$
are pairwise distinct.  Facts~\ref{fact: min length} and~\ref{fact:
distinct prefixes} show that the proportion of elements of $R_n^k$ in
$Y_{\alpha,\lambda,n,k}$ converges to 1 exponentially fast: in the
search for exponentially generic properties of subgroups, we can
restrict our attention to the tuples in $Y_{\alpha,\lambda,n,k}$ and
to the subgroups they generate.

The following observation is the basis for our exponential genericity 
proofs, in the context of the word-based distribution.

\begin{fact}\label{fact: essential word-based}
    Let $\alpha, \lambda$ satisfy $0<2\lambda < \alpha < 1$.  If $\vec
    h\in Y_{\alpha,\lambda,n,k}$ and $H = \langle\vec h\rangle$, then
    $\Gamma(H)$ consists of two parts:
    
    - the vertices at distance at most $\lceil\lambda n\rceil$ from
    the distinguished vertex and the edges connecting them: this forms
    a tree with vertex set $\Pref_{\lambda n}(\vec h)$, and edges $u
    \to^a ua$ if $a\in A$ and $u,ua\in \Pref_{\lambda n}(\vec h)$;
    this tree, which we call the \textit{central part} of $\Gamma(H)$,
    has $2k$ leaves;
    
    - and for each $1\le i\le k$, where $h_i = p_im_is_i$ and $|p_i| =
    |s_i| = \lceil\lambda n\rceil$, a path labeled $m_i$ from the
    vertex $p_i$ to the vertex $s_i\inv$ (both are in the central
    part); we call these paths the \textit{outer loops}.
\end{fact}

This leads to the following results.  Propositions~\ref{prop rank k}
and~\ref{distinct subgroups} first appeared in a paper by Jitsukawa
\cite{2002Jitsukawa}.  They are direct consequences of earlier
analogous results (counting cyclic words instead of words) due to
Arzhantseva and Ol'shanski\u\i\ \cite[Lemma 3]{1996ArjantsevaO}.

\begin{proposition}\label{prop rank k}
    Exponentially generically, a $k$-tuple of elements of $R_n$
    generates a subgroup of rank $k$.
\end{proposition}

\proof
Let $\alpha, \lambda$ satisfy $0<2\lambda < \alpha < 1$.  As observed
above, it suffices to show that if $\vec h\in Y_{\alpha,\lambda,n,k}$,
then $H = \langle\vec h\rangle$ has rank $k$.  In that case, using
Fact~\ref{fact: essential word-based}, we find that $\Gamma(H)$ is
formed of a central part and $k$ outer loops.

The central part is a tree and like all trees, the number of its edges
is 1 less than the number of its vertices.  With the notation of
Fact~\ref{fact: essential word-based}, the number of additional
vertices (resp.  edges) in the outer loops is $\sum_i(|m_i|-1)$ (resp.
$\sum_i|m_i|$).  Therefore, in $\Gamma(H)$, we have $|E|-|V|+1 = k$ and hence $\rank(H) = k$ (see Section~\ref{sec: stallings}).
\eop

\begin{proposition}\label{corollary free product}
    Exponentially generically, a $k$-tuple $\vec h$ and an
    $k'$-tuple $\vec h'$ of elements of $R_n$ generate subgroups
    that are distinct, have trivial intersection, and are such that
    $\langle \vec h, \vec h'\rangle = \langle\vec h\rangle \ast
    \langle\vec h'\rangle$.
\end{proposition}

\proof
Since the first $k$ components of a $(k+k')$-tuple of elements of
$R_n$ are independent from the $k'$ last components, and since such
a $(k+k')$-tuple exponentially generically generates a subgroup of
rank $k+k'$ (Proposition~\ref{prop rank k}), we find that a
$k$-tuple and an $k'$-tuple of elements of $R_n$ exponentially
generically generate their free product.  This in turn implies the
other properties.
\eop

Proposition \ref{corollary free product} shows that two $k$-tuples of
elements of $R_n$ exponentially generically generate distinct
subgroups.  Proposition~\ref{distinct subgroups} is a little more
precise.

\begin{proposition}\label{distinct subgroups}
    Let $\alpha, \lambda$ satisfy $0<2\lambda < \alpha < 1$.  The
    $k$-tuples $\vec h$ and $\vec h'$ in $Y_{\alpha,\lambda,n,k}$
    generate distinct subgroups, unless $\vec h' =
    (h_{\sigma(1)}^{\epsilon_1}, \ldots, h_{\sigma(k)}^{\epsilon_k})$
    for some permutation $\sigma$ of $[k]$ and for
    $\epsilon_1,\ldots,\epsilon_k \in \{-1,+1\}$.
\end{proposition}

\proof
If $\langle\vec h\rangle = \langle\vec h'\rangle$, then the graphs
$\Gamma\langle\vec h\rangle$ and $\Gamma\langle\vec h'\rangle$ are
equal.  In particular, their central parts, formed by the vertices at
distance at most $\lceil\lambda n\rceil$ from the distinguished
vertex, coincide.  Fact~\ref{fact: essential word-based} shows that if $\vec h \in Y_{\alpha,\lambda,n,k}$, then $\Gamma\langle\vec h\rangle$ completely determines all the $k$-tuples $\vec h' \in Y_{\alpha,\lambda,n,k}$ such that $\Gamma\langle\vec h'\rangle = \Gamma\langle\vec h\rangle$, and that they coincide with $\vec h$ up to the order of the elements and the direction in which the outer loops are read.
\eop

Proposition~\ref{distinct subgroups} shows that, if we consider the
class $S_{\alpha,\lambda,n,k}$ of subgroups generated by $k$-tuples in
$Y_{\alpha,\lambda,n,k}$, then each subgroup occurs the same number of
times, namely $2^kk!$.  Randomly choosing a $k$-tuple in
$Y_{\alpha,\lambda,n,k}$ yields therefore a random subgroup in
$S_{\alpha,\lambda,n,k}$, and the proportion of these subgroups among
all subgroups generated by a $k$-tuple of words of length at most $n$
tends to $1$ exponentially fast.

\subsection{The graph-based distribution}\label{sec: graph based}

The uniform distribution on the set of size $n$ Stallings graphs was
analyzed by Bassino, Nicaud and Weil \cite{BNW}.  Here we summarize
the principles of this distribution and the features which
will be used in this paper.

We already noted that in Stallings graphs, each letter labels a
partial injection on the vertex set: in fact, a Stallings graph can be
viewed as a collection $(f_a)_{a\in A}$ of partial injections on an
$n$-element set, with a distinguished vertex, and such that the
resulting graph (with an $a$-labeled edge from $i$ to $j$ if and only
if $j = f_a(i)$) is connected and has no vertex of degree 1, except
perhaps the distinguished vertex.  We may even assume that the
$n$-element set in question is $[n] = \{1,\ldots,n\}$, with 1 as the
distinguished vertex, see \cite[Section 1.2]{BNW} for a precise
justification. In particular, the automorphism group of an admissible pair $(\Gamma,1)$ is always trivial.

One shows \cite[Corollary 2.7]{BNW} that the probability that an
$A$-tuple $(f_a)_{a\in A}$ of partial injections on $[n]$ induces a
Stallings graph tends to 1 as $n$ tends to infinity, and the problem
of randomly generating a Stallings graph then reduces (via an
efficient rejection algorithm, see \cite[Section 3]{BNW}) to the
problem of efficiently generating a random partial injection on $[n]$.
This view of a Stallings graph as an $A$-tuple of partial injections
on $[n]$ is central in our analysis.

The maximal orbits of a partial injection $f$ (equivalently: the
connected components of the function graph of $f$) can be of two
kinds: cycles -- where each element is both in the domain and in the
range of $f$ -- and sequences.  The size of each of these components
is defined to be the number of vertices which they contain.  It is
this combinatorial view of partial injections -- as a disjoint union
of cycles and sequences --, which is at the heart of the random
generation algorithm, obtained using the so-called recursive method
\cite{NijenhuisWilf78,FlajoletZvC94}.

The distribution of sizes of components is studied in \cite[Section
3]{BNW}, as well as the distribution of cycles vs.  sequences among
size $k$ components.  The random generation algorithm consists in
drawing a size of component, say $k$, according to the relevant
distribution; then drawing whether this size $k$ component is a cycle
or a sequence; and finally drawing a partial injection on the
remaining $n-k$ elements \cite[Section 3.1]{BNW}.  This results in a
partial injection on an $n$-element set, and we need only add a random
numbering (1 through $n$) of the elements of that set.

However complex the method may seem, it guarantees a uniform
distribution among all size $n$ partial injections, it is easy to
implement and its average time complexity is linear (in the RAM model;
it is $\O(n^2\log n)$ under the bit-cost assumption) \cite[Section
3.3]{BNW}.

To further discuss partial injections and other combinatorial
structures, we use the notion of \textit{exponential generating
series}, written EGS. If $\calS$ is a class of finite discrete
structures such that there are finitely many $\calS$-structures of
each size, let $S_n$ be the number of $\calS$-structures of size $n$.
The EGS of $\calS$ is the formal power series $S(z) = \sum_{n\ge
0}\frac{S_n}{n!}z^n$.

Let $I(z) = \sum_{n\ge 0}\frac{I_n}{n!}z^n$ be the EGS of partial
injections.  Bassino, Nicaud and Weil show the following \cite[Section
2.1 and Proposition 2.10]{BNW}.

\begin{proposition}\label{prop: facts graphs}
    The EGS $I(z)$ of partial injections satisfies the following
    $$I(z) = \frac1{1-z}\exp\left(\frac z{1-z}\right) \quad
    \textrm{and} \quad \frac{I_n}{n!} =
    \frac{e^{-\frac12}}{2\sqrt\pi}e^{2\sqrt
    n}n^{-\frac14} (1 + o(1)).$$
\end{proposition}

This result is obtained by means of deep theorems from analytic
combinatorics.  The same methods can be used to study the asymptotic
behavior of particular parameters, such as the number of sequences of
a partial injection.  This parameter is directly connected with the
number of edges in the Stallings graph formed by the partial
injections $f_a$ ($a\in A$), which leads to the following result
\cite[Lemma 2.11 and Corollary 4.1]{BNW}.

\begin{proposition}\label{prop: expected sequences}
    The expected number of sequences in a randomly chosen partial
    injection of size $n$ is asymptotically equivalent to $\sqrt n$.
    
    The expected rank of a randomly chosen size $n$ subgroup of $F_r$
    is asymptotically equivalent to $(r-1)n-r\sqrt n+1$.
\end{proposition}

\subsection{Negligible and generic properties of subgroups}
\label{sec: begin discussion}

Thus, in the discussion of statistical properties of finitely
generated subgroups of a (fixed) free group $F_r$, we have two
distributions at our disposal.  One, the word-based distribution, is
governed by two parameters: the number of generators and their
maximum length, the former fixed and the latter allowed to tend to
infinity. The other, the graph-based distribution, is governed by a
single parameter: the size of the Stallings graph.

We first observe that our discussion of the graph-based distribution
(as well as the results in \cite{BNW}) is in terms of spheres rather
than balls: as we saw in the Section~\ref{sec: balls spheres}, the
(exponential) negligibility or genericity results obtained in that
setting are sufficient.  In contrast, the existing literature on the
word-based distribution is in terms of balls, as is our description in
Section~\ref{sec: word based} above.

The graph-based as well as the word-based distribution allow the
discussion of properties of subgroups (of subgroups of a fixed rank
$k$ in the word-based case).  There is of course no reason why a
property that is generic or negligible in one distribution should have
the same frequency in the other.

\begin{figure}
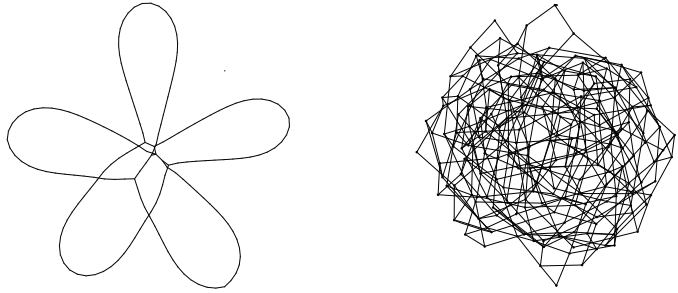

  \begin{center}
   \hspace{1.5cm}\begin{minipage}[c]{.4\linewidth}
     \rotatebox{0}{ \includegraphics[scale=0.5]{jits.ps}}
   \end{minipage}
   \begin{minipage}[c]{.4\linewidth}
     \rotatebox{0}{ \includegraphics[scale=0.5]{example.ps}}
   \end{minipage}
   \caption{\small The Stallings graphs of two randomly generated
   subgroups of $F_2$.  On the left, a subgroup generated by a random
   5-tuple of words of length at most 40.  On the right, a random
   Stallings graph of size 200.  Only the shape of the graphs is
   depicted, vertices and edge labels and directions are not
   represented.  The pictures have been generated by {\tt neato}.
   Note that the scale (average distance between two vertices) is not
   the same on the two pictures.\label{figure 200}}
 \end{center}
\end{figure}

Our two distributions are indeed very different.  How different is
illustrated in Figure~\ref{figure 200}, which shows a ``random'' size
200 Stallings graph and the Stallings graph of the subgroup of $F_2$
generated by a ``random'' 5-tuple of words of length at most 40 (which
has close to 200 vertices).  This figure provides the intuition to
exhibit properties of subgroups that are negligible in one
distribution and generic in the other.

It is not difficult to come up with such properties.  It is the case,
for instance, of the property \textit{to have rank $\ell$}, for a
fixed integer $\ell\ge 1$.  In the graph-based distribution, this
property is negligible as a consequence of Proposition~\ref{prop:
expected sequences} (see \cite[Corollary 4.2]{BNW}).  In contrast, it
is exponentially generic in the word-based distribution \textit{of
$\ell$-generated} subgroups, see Proposition~\ref{prop rank k}.  For
the same reason, it is exponentially negligible in the word-based
distribution \textit{of $k$-generated} subgroups with $k \ne \ell$.

The properties of malnormality and purity, discussed in
Section~\ref{sec: malnormal}, provide more complex examples of this
sort.

\section{Malnormal and pure subgroups}\label{sec: malnormal}

Malnormality and purity are two important properties of subgroups.  A
subgroup $H$ is \textit{pure} if $x^n\in H$ and $n\ne 0$ implies $x\in
H$.  A pure subgroup is also called \textit{closed under radicals} or 
\textit{isolated}.

The subgroup $H$ is \textit{malnormal} if $H \cap H^g = 1$ for every
$g\not\in H$.  Malnormal subgroups play an important role in the study
of amalgamated products (e.g.
\cite{1971KarrassSolitar,1998BleilerJones}) and in the characterization
of their hyperbolicity \cite{1998KharlampovichMyasnikovTAMS2}.  The
following is elementary from the definition.

\begin{lemma}\label{malnormal 2 pure}
    If a subgroup is malnormal, then it is pure.
\end{lemma}

Note that the converse statement does not hold: $\langle
a,bab\inv\rangle$ is pure, yet not malnormal.

Both malnormality and purity have nice graphical characterizations,
which imply that these properties are decidable for finitely generated
subgroups of free groups.  The result on malnormality is due to
Kapovich and Myasnikov \cite{2002KapovichMyasnikov} (following a
decidability result in \cite{1999Baumslag}), that on purity is due to
Birget, Margolis, Meakin and Weil \cite{2000BirgetMMW}.

\begin{proposition}\label{prop: charact malnormal}
    Let $(\Gamma,1)$ be the graphical representation of a subgroup $H$.
    
    \begin{itemize}
	\item[(1)] $H$ is non-malnormal if and only if there exists a
	non-trival reduced word $u$ and distinct vertices $x\ne y$ in
	$\Gamma$ such that $u$ labels loops at $x$ and at $y$.
	
	\item[(2)] $H$ is non-pure if and only if there exists a non-trival
	reduced word $u$, an integer $n\ge 2$ and a vertex $x$ in
	$\Gamma$ such that $u^n$ labels a loop at $x$ but $u$ does
	not.
    \end{itemize}
\end{proposition}

\subsection{Genericity in the word-based distribution\dots}
\label{generic malnormal}

Jitsukawa shows that malnormality is a generic property in free 
groups \cite[Theorem 4 and Lemma 6]{2002Jitsukawa}. His arguments can be 
extended to show that it is exponentially generic.

\begin{theorem}\label{thm: malnormality word-based}
    Malnormality is exponentially generic in the word-based distribution.
\end{theorem}

In view of Lemma~\ref{malnormal 2 pure}, we also have the following result.

\begin{corollary}
    Purity is exponentially generic in the word-based distribution.
\end{corollary}

We now proceed to prove Theorem~\ref{thm: malnormality word-based}.
The proof relies on the two following lemmas, which provide an
analogue of a small cancellation property for tuples of reduced words.
These lemmas constitute a variant of the results of Gromov \cite[0.2.A]{1987Gromov} for tuples of cyclic words, proved also by Champetier \cite{1991Champetier,1995Champetier} and by Arzhantseva
and Ol'shanski\u\i\ \cite[Lemma
3]{1996ArjantsevaO}.

\begin{lemma}\label{fact: C' single word}
    Let $0 < \beta < 1$.  The proportion of $k$-tuples $\vec h$ of
    reduced words in $F_r$ of length at most $n$, such that one of the
    $h_i$ contains two distinct occurrences of factors $v$ and $w$ of
    length at least $\beta n$, with $v = w$ or $v = w\inv$, converges
    to 0 exponentially fast.
\end{lemma}

\begin{lemma}\label{fact: C' double word}
    Let $0 < \beta < 1$.  The proportion of $k$-tuples $\vec h$ of
    reduced words of length at most $n$, such that a word $v$ of
    length at least $\beta n$ has an occurrence in one of the $h_i$
    and $v$ or $v\inv$ has an occurrence in $h_j$ for some $j\ne i$,
    converges to 0 exponentially fast.
\end{lemma}


We can now prove Theorem~\ref{thm: malnormality word-based}.

\proofof{Theorem~\ref{thm: malnormality word-based}}
Let $0 < \alpha < 1$ and $0 < \lambda < \frac\alpha4$.  Let $\vec h =
(h_1,\ldots,h_k)$.  Exponentially generically, we have $\min|h_i| >
\alpha n$ and the prefixes of length $\lceil\lambda n\rceil$ of the
$h_i$ and $h_i\inv$ are pairwise distinct (Facts~\ref{fact: min
length} and~\ref{fact: distinct prefixes}).  In addition,
exponentially generically, no word of length at least
$\frac{\alpha-4\lambda}2n$ has distinct occurrences as a factor of the
$h_i$ and the $h_i\inv$ (Lemmas~\ref{fact: C' single word}
and~\ref{fact: C' double word}).

Let us now assume that $\vec h$ satisfies all these properties.  Then
$\Gamma = \Gamma(\langle\vec h\rangle)$ is composed of a central part,
which is a tree containing the distinguished vertex and all the
vertices corresponding to the prefixes of the $h_i$ and $h_i\inv$ of
length up to $\lceil\lambda n\rceil$, and of outer loops whose labels
are factors of the $h_i$ (or the $h_i\inv$, depending on the direction
in which they are read), see Fact~\ref{fact: essential word-based}.

Any loop in $\Gamma$ must visit the central part of $\Gamma$ at least
once, and run along at least one of the outer loops.  Let us now
assume that a word $u$ labels two distinct loops in $\Gamma$.  Up to
conjugation of $u$, we can assume that the base point of the first
loop is in the central part of $\Gamma$.  Then $u$ has a factor $v$ of
length $\lceil\alpha n\rceil - 2\lceil\lambda n\rceil$, which is a
factor of some $h_i$ or $h_i\inv$.  The other occurrence of a loop
labeled $u$ reveals another path in $\Gamma$ labeled $v$.  This path
may not be entirely in an outer loop, but if it is not, then it visits
the central part of $\Gamma$ only once, so it has a factor $v'$ of
length at least $\frac{\alpha-4\lambda}2n$ in an outer loop, and hence
in one of the $h_i$ or $h_i\inv$.  This word $v'$ has distinct
occurrences as a factor of the $h_i$ and the $h_i\inv$, a
contradiction.
\eopo

\subsection{\dots\ and negligibility in the graph-based distribution}\label{sec: malnormal graph}

In contrast, we show that malnormality and purity are negligible in
the graph-based distribution.

\begin{theorem}\label{thm: purity negligible}
    The probability that a random subgroup of size $n$ is pure is
    $\O(n^{-\frac r2})$.
\end{theorem}

By Lemma~\ref{malnormal 2 pure}, this implies the following

\begin{corollary}\label{thm: malnormality negligible}
    The probability that a random subgroup of size $n$ is malnormal is
    $\O(n^{-\frac r2})$.
\end{corollary}

To prove Theorem~\ref{thm: purity negligible}, we observe that if $H$
is a finitely generated subgroup of $F_r$ and a
cycle of length at least 2 in $\Gamma(H)$ is labeled by a power of some letter $a$, then $H$ is not pure
(Proposition~\ref{prop: charact malnormal}).  Therefore, if a subgroup
is pure, then the partial injection determined by each letter in $A$
has only sequences and length 1 cycles.

Thus Theorem~\ref{thm: purity negligible} follows directly from the 
following proposition.

\begin{proposition}\label{propK}
    The probability that a size $n$ partial injection has no cycle of
    length greater than or equal to 2 is asymptotically equivalent to
    $\frac e{\sqrt n}$.
\end{proposition}

Our proof of Proposition~\ref{propK} uses Hayman's theorem, discussed 
in Section~\ref{sec: Hayman} below.

\begin{remark}\label{remark purity exp negligible}
    There are many more reasons for a subgroup to fail to be pure,
    than those considered here.  In terms of Proposition~\ref{prop:
    charact malnormal}, we have considered only the words $u$ that are
    equal to a letter of the alphabet.  As a result, the probability
    of purity and that of malnormality are likely to be much smaller
    than the upper bounds given above.  The open question here is
    whether purity and normality are exponentially negligible with
    respect to the graph-based distribution.
\end{remark}

\subsubsection{H-admissible functions and Hayman's theorem}
\label{sec: Hayman}

Hayman's theorem on the asymptotic behavior of the coefficients of
certain power series requires a technical hypothesis called
H-admissibility.  Here we give only the technical definition and
statement we will use, and we refer the readers to \cite[Chapter
VIII]{FlajoletSedgewick} for further details on this theorem and on
saddlepoint asymptotics in general.

Let $f(z)$ be a function of the form $f(z) = e^{h(z)}$ that is
analytic at the origin, with radius of convergence $\rho$.  We denote
by $[z^n]f(z)$ the coefficient of $z^n$ in the power series
development of $f$ at the origin.  Let
$$a(r)=rh'(r) \quad \mbox{and} \quad b(r)=r^2h^{''}(r)+rh'(r).$$
The function $f(z)$ is said to be {\it H-admissible} if there exists a
function $\delta\colon ]0,\rho[ \longrightarrow ]0,\pi[$ such that the
following three conditions hold:
\begin{itemize}
    \item[(H1)] $\lim_{r \rightarrow \rho} b(r) = +\infty.$
    \item[(H2)] Uniformly for $|\theta| \leq \delta(r)$
    $$f(re^{i\theta}) \sim f(r)e^{i\theta a(r)
    -\frac{1}{2} \theta^2 b(r)}\quad\hbox{when $r$ tends to $\rho$.}$$
    [That is, $f(re^{i\theta}) = f(r)e^{i\theta a(r)
    -\frac{1}{2} \theta^2 b(r)} (1 + \gamma(r,\theta))$ with 
    $|\gamma(r,\theta)| \le \tilde\gamma(r)$ when $|\theta| \le 
    \delta(r)$ and $\lim_{r \rightarrow \rho}\tilde\gamma(r) = 0$.]
    \item[(H3)] and uniformly for $\delta(r)\leq |\theta| \leq \pi$ 
    $$f(re^{i\theta})\sqrt{b(r)} = o(f(r)) \qquad \mbox{when $r$ 
    tends to $\rho$.}$$
\end{itemize}

Hayman's theorem \cite[Theorem VIII.4]{FlajoletSedgewick} states the following.

\begin{theorem} \label{thcol}
    Let $f(z) = e^{h(z)}$ be a H-admissible function with radius of 
    convergence $\rho$ and
    $\zeta=\zeta(n)$ be the unique solution in the interval $]0,\rho[$
    of the saddlepoint equation
    $$\zeta \frac{f'(\zeta)}{f(\zeta)}=n.$$
    Then
    $$[z^n]f(z) = \frac{f(\zeta)}{\zeta^n \, \sqrt{2 \pi
    b(\zeta)}}\left(1+o(1)\right).$$
    where $b(z) = z^2h''(z)+zh'(z)$.
\end{theorem}

\subsubsection{Proof of Proposition~\ref{propK}}\label{proofofpropK}

Let $\calK$ be the set of partial injections in which all the cycles
have length 1 and let $\calJ$ the set of partial injections without
any cycles (a subset of $\calK$).  The elements of $\calJ$ are known
as \textit{fragmented permutations}, see \cite[Section
II.4.2]{FlajoletSedgewick}.

Let $K_n$ and $J_n$ be the number of size $n$ elements of $\calK$ and
$\calJ$, and let $K(z)$ and $J(z)$ be the corresponding EGS. We want to show that $\frac{K_n}{I_n}$ is equivalent to $\frac e{\sqrt n}$.

The
series $J(z)$ is studied in detail in \cite[Example VIII.7, Proposition
VIII.4]{FlajoletSedgewick}.  There, it is shown in particular that
$J(z)$ is H-admissible and that
\begin{equation}
    J(z) = \exp\left(\frac z{1-z}\right) \quad\textrm{and}\quad\frac{J_n}{n!} =
    \frac{e^{-\frac12}}{2\sqrt\pi}e^{2\sqrt n}n^{-\frac34} (1 + 
    o(1)).
\end{equation}

A partial injection in $\calK$ consists of a set of length 1 cycles 
and a fragmented permutation. It follows that
$$K_n = \sum_{k=0}^n \frac{n!}{k!(n-k)!}J_k,$$
so that
$$K(z) = \sum_{n=0}^\infty\frac{K_n}{n!}z^n =
\left(\sum_{n=0}^\infty\frac1{n!}z^n\right)
\left(\sum_{n=0}^\infty\frac{J_n}{n!}z^n\right)
= e^zJ(z) = \exp\left(z+\frac{z}{1-z}\right).$$
Now $e^z$ is H-admissible: this can be verified directly, or by
application of \cite[Theorem VIII.5]{FlajoletSedgewick}.  We already
noted that $J(z)$ is H-admissible, and hence $K(z)$ is H-admissible as
well, as the product of two H-admissible functions (\cite[Theorem
VIII.5]{FlajoletSedgewick} again).

The saddle-point equation $z\frac{K'(z)}{K(z)} = n$ (see
Section~\ref{sec: Hayman}) is
\begin{align*}
    \frac{z(2-2z+z^2)}{(1-z)^2} &= n, \\
    \textit{i.e.}\qquad z^3 - (n+2) z^2 + 2(n+1)z - n &= 0.
\end{align*}
Let $P_n(z)$ be the polynomial on the left hand side of this last
equation.  Examining the sign of the derivative of $P_n(z)$ on the
interval $[0,1]$ and the values of $P_n$ at $0$ and $1$, we find that
$P_n$ has a unique zero between $0$ and $1$, say $\zeta_n$.  Moreover
\begin{equation}\label{zeta development}
    \zeta_n = 1 - \frac1{\sqrt n} +\frac1{2n}+\O\left(\frac1{n\sqrt n}\right).
\end{equation}
This asymptotic development can be obtained using \textsc{maple}, 
based on the application of the Cardan method to this degree $3$ 
polynomial. We can also observe the following. Let $Q_n(z)$ be the 
polynomial defined by the identity
$$P_n(1-z) = 1 - z + (1 - n) z^2 - z^3 = Q_n(z) - z^3.$$
The zero of $Q_n(z)$ in the interval $[0,1]$ is
$$\alpha_n = \frac{\sqrt{4n-3}-1}{2(n-1)} = \frac1{\sqrt n} - 
\frac1{2n} + \O(n^{-\frac32})$$
and if $\beta_n = 1-\alpha_n$, we have $P_n(\beta_n) = -\alpha_n^3$,
which is negative for $n$ large enough.

Now let $\gamma_n = 1 - \frac1{\sqrt n} + \frac1{2n}$. Then
$$P_n(\gamma_n) = \frac5{4n} + \O(n^{-\frac32}),$$
which is positive for $n$ large enough.  It follows that $\beta_n <
\zeta_n < \gamma_n$, justifying the development in (\ref{zeta
development}).

With the notation of Section~\ref{sec: Hayman}, we also have
\begin{align*}
    b(z) &= z^2\frac{d^2}{dz^2}\log K(z) + z\frac{d}{dz}\log K(z)\\
    &= \frac{z(2-2z+3z^2-z^3)}{(1-z)^3},\\
    \textrm{so that }b(\zeta_n) &= 2n^{3/2} +\O(n).
\end{align*}
Then we have
\begin{align*}
    \log \zeta_n^n &= n\log\left(1 - \frac1{\sqrt n}
    +\frac1{2n}+\O\left(\frac1{n\sqrt n}\right)\right) = -\sqrt{n}
    + \O\left(\frac1{\sqrt n}\right)\\
    \textrm{and }K(\zeta_n) &= \exp\left(\zeta_n + \frac{\zeta_n}{1-\zeta_n}\right)\sim 
    e^{\frac12}e^{\sqrt n}.
\end{align*}
By Theorem~\ref{thcol}, we now have
$$[z^n]K(z) \sim \frac{K(\zeta_n)}{\zeta_n^n \sqrt{2\pi b(\zeta_n)}}
\sim e^{\sqrt n}e^{\frac12}e^{\sqrt
n}\frac1{2\sqrt{\pi}n^{3/4}}
\sim \frac{e^\frac12}{2\sqrt{\pi}}n^{-3/4}e^{2\sqrt n}.$$
Proposition~\ref{propK} follows since $[z^n]I(z)\sim \frac{e^{-\frac12}}{2\sqrt{\pi}}e^{2\sqrt n}n^{-1/4}$ (Proposition~\ref{prop: facts graphs}) and hence, $\frac{K_n}{I_n} \sim \frac e{\sqrt n}$.

\subsection{A remark on the Hanna Neumann conjecture}\label{sec: SHNC}

The Hanna Neumann Conjecture (HNC), recently established by Mineyev \cite{2011-Mineyev-Annals} after several decades of partial results, deals with the rank of the
intersection of finitely generated subgroups of free groups (see also \cite{2011-Mineyev-preprint} for an alternative proof, purely in terms of groups and graphs). For convenience, let the
\textit{reduced rank} of a subgroup $H$, written $\rrank(H)$, be equal
to
$$\rrank(H) = \max(0, \rank(H)-1).$$
Mineyev's theorem states that, if $H$ and $K$ are finitely generated subgroups 
of $F$, then $\rrank(H\cap K) \le
\rrank(H)\rrank(K)$, as conjectured by Hanna Neumann. It also shows the stronger inequality conjectured by Burns \cite{1971Burns}, formerly known as the strengthened Hanna Neumann conjecture (SHNC):
$$\sum \rrank(H \cap K^g) \le \rrank(H)\rrank(K),$$
where the sum runs over all subsets $KgH$ (in $K\backslash F/H$) such that $H \cap
K^g \ne 1$ and $K^g = g\inv K g$.

We discuss here how our results show that the cases where this inequality is non-trivial are rare, in the sense that $\rrank(H \cap K)$ and $\sum \rrank(H \cap K^g)$ are generically equal to zero. 

It was observed, initially by Stallings \cite{1983Stallings} and
Gersten \cite{1983Gersten}, that HNC and SHNC have natural
interpretations in terms of Stallings graphs. If $\Gamma$ is an 
$A$-labeled graph, let us denote by $\chi(\Gamma)$ the 
difference between the number of edges and the number of vertices of 
$\Gamma$: thus $\rrank(H) = \chi(\Gamma(H))$.  Let $\Delta(H,K)$ be
the graph obtained from $\Gamma(H)$ and $\Gamma(K)$ as follows: the
vertices of $\Delta(H,K)$ are the pairs $(u,v)$ such that $u$ is a vertex
of $\Gamma(H)$ and $v$ is a vertex of $\Gamma(K)$; and the edges of
$\Delta(H,K)$ are the triples $((u,v),a,(u',v'))$ such that $(u,a,u')$ is an
edge of $\Gamma(H)$ and $(v,a,v')$ is an edge of $\Gamma(K)$.

Let $\Delta_1$ be the connected component of $\Delta(H,K)$ containing
$(1,1)$ (where 1 denotes the origin of $\Gamma(H)$ and of
$\Gamma(K)$), and let $\Delta_2$ be the union of the connected
components of $\Delta(H,K)$ which are not trees.  Then HNC holds for
$H$ and $K$ if and only if $\chi(\Delta_1) \le \rrank(H)\rrank(K)$,
and SHNC holds for $H$ and $K$ if and only if $\chi(\Delta_2) \le
\rrank(H)\rrank(K)$.

Now observe (as in Proposition~\ref{corollary free product}) that a 
randomly chosen $(k+\ell)$-tuple of elements of $R_n$ is composed of 
the juxtaposition of a randomly chosen $k$-tuple and a randomly 
chosen $\ell$-tuple. Exponentially generically, such a $k$-tuple 
$\vec h$ and $\ell$-tuple $\vec h'$ generate subgroups with trivial 
intersection: in particular, HNC holds exponentially generically in 
the word-based distribution.

In fact, with the same ideas as in the proof of Theorem~\ref{thm:
malnormality word-based}, exponentially generically, there is no loop
in $\Gamma(\langle\vec h\rangle)$ with an occurrence as a loop in
$\Gamma(\langle\vec h'\rangle)$.  Therefore SHNC holds exponentially
generically.

As mentioned above, both HNC and SHNC are now known to hold, but it seems interesting to point out that exponentially generically, they hold for trivial reasons.

\section{An intermediate property}\label{sec: intermediate}

In this section, we discuss an \textit{intermediate} property of
subgroups, that is a property such that the proportion of subgroups of
size $n$ with this property has a limit which is neither 0 nor 1
(respectively the negligible and the generic cases).

\begin{theorem}\label{thm: visible property}
    The probability that a random size $n$ subgroup of $F_r$ intersects 
    trivially the conjugacy classes of the generators tends to 
    $e^{-r}$ when $n$ tends to infinity.
\end{theorem}

The discussion of this property is included here because we do not
know many examples of such intermediate properties.  Unfortunately, the property
in question is geometric in the sense that it depends on the
combinatorial parameters of the Stallings graph of the subgroup, and
is not preserved under the automorphisms of $F_r$.  It would be
interesting to exhibit such a property that would be algebraic
(preserved under automorphisms).  One might think for instance of the
property of avoiding the conjugacy classes of all the elements of some
basis of $F_r$, or the property of avoiding all primitive words.

\begin{remark}
    The property described in Theorem~\ref{thm: visible property} is
    exponentially negligible in the word-based distribution.  Indeed, 
    if $\vec h$ is a $k$-tuple of reduced words of length at most 
    $n$, then
    $\Gamma(\langle\vec h\rangle)$ has exponentially generically $k$ 
    loops of length at least $\frac n2$ and no loop of length 1 (see 
    the discussion in Section~\ref{sec: word based} with $\alpha = 
    \frac34$ and $\lambda = \frac18$).
\end{remark}

The rest of this section is devoted to the proof of Theorem~\ref{thm:
visible property}.  It is easily verified that a subgroup $H$ contains
a conjugate of letter $a\in A$ if and only if $a$ labels a loop at
some vertex of $\Gamma(H)$, that is, if and only if the corresponding
partial injection has some fixpoint.  Since the drawing of the partial
injections corresponding to the different letters is independent, the
theorem follows directly from the following proposition.

\begin{proposition}\label{prop: no fixpoint}
    The probability that a size $n$ partial injection has no fixpoint
    tends to $\frac1e$ when $n$ tends to infinity.
\end{proposition}

\begin{remark}
    Note that $\frac1e$ is also the limit of the probability that a
    size $n$ permutation has no fixpoint (a so-called
    \textit{derangement}, see \cite{Comtet1974}).
\end{remark}

Our proof of Proposition~\ref{prop: no fixpoint} again uses Hayman's
theorem (Section~\ref{sec: Hayman}).  We also need the following
technical result.

\begin{proposition}\label{lemma: no fixpoint}
    Let $f_0(z)$ be an H-admissible function with radius of
    convergence $\rho < \infty$.  Then $f(z) = e^{-z}f_0(z)$ is
    H-admissible as well.
\end{proposition}

\proof
Since $f_0$ is analytic at the origin, it is clear that $f(z)$ is
analytic at the origin as well, with a radius of convergence equal to
that of $f_0(z)$.

Let $h(z)$ be such that $f(z) = e^{h(z)}$.  If $h(z) = h_0(z) - z$,
then we have $f_0(z) = e^{h_0(z)}$.

Let $a_0(t)  = th_0'(t)$, $a(t) = th'(t)$, $b_0(t) = t^2h''_0(t) + 
a_0(t)$ and $b(t) = t^2h''(t) +a(t)$. Then $a(t) = a_0(t) - t$ and 
$b(t) = b_0(t) - t$.

It is immediate that $\lim_{t \to \rho}b(t) = +\infty$ since this
limit holds for $b_0$.  That is, Condition (H1) holds.

We now verify Condition (H2).  Let $\delta(t)$ be a positive function
such that $\lim_{t \to \rho}\delta(t) = 0$; and such that, uniformly
for $|\theta| \le \delta(t)$, and as $t$ tends to $\rho$,
$$h_0(te^{i\theta}) = h_0(t) + i\theta a_0(t) - \frac12\theta^2 b_0(t)
+o(1).$$
Then
\begin{align*}
    h(te^{i\theta}) &= h_0(te^{i\theta}) - te^{i\theta} \\
    &= h_0(t) + i\theta a_0(t) - \frac12\theta^2 b_0(t) +o(1) -
    te^{i\theta}\\
    &= h(t) + i\theta a(t) - \frac12\theta^2 b(t) +o(1) -
    te^{i\theta} + t + ti\theta - \frac12t\theta^2.
\end{align*}
We now observe that, if $|\theta|\le \delta(t)$ and as $t$ tends to
$\rho$, then $|ti\theta| \le t\delta(t) = o(1)$ and similarly,
$\frac12t\theta^2 = o(1)$.  Finally,
$$|t(1-e^{i\theta})| = t \sqrt{(1- \cos\theta)^2+\sin^2\theta} =
t\sqrt{2(1- \cos\theta)} \le t|\theta| \le t\delta(t) = o(1).$$
Thus $h(te^{i\theta}) = h(t) + i\theta a(t) - \frac12\theta^2 b(t) + 
o(1)$ uniformly for $|\theta| \le \delta(t)$, which concludes the verification of (H2).

Finally, we want to show that
$\frac{f(te^{i\theta})\sqrt{b(t)}}{f(t)}$ tends to $0$ when $t$ tends
to $\rho$, uniformly for $\delta(t) \le |\theta| \le \pi$. We have
\begin{align*}
    \frac{f(te^{i\theta})\sqrt{b(t)}}{f(t)} &=
    \frac{f_0(te^{i\theta})e^{-te^{i\theta}}\sqrt{b(t)}}{f_0(t)e^{-t}} \\
    &= \frac{f_0(te^{i\theta})\sqrt{b_0(t)}}{f_0(t)} e^{t(1- e^{i\theta})}
    \sqrt{1 - \frac t{b_0(t)}}.
\end{align*}
Since $f_0$ is H-admissible, uniformly for $\delta(t) \le |\theta| 
\le \pi$ and as $t$ tends to $\rho$,
$$\frac{f_0(te^{i\theta})\sqrt{b_0(t)}}{f_0(t)} = o(1).$$
Moreover, $\sqrt{1 - \frac{t}{b_0(t)}} = 1 + o(1)$ since $\lim_{t\to
\rho}b_0(t) = +\infty$.  Finally, when $0 < t < \rho$,
$|e^{t(1-e^{i\theta})}| = e^{t(1-\cos\theta)} \le e^{2\rho}$.  This
suffices to conclude that (H3) holds, and hence that $f(z)$ is
H-admissible.
\eop

\proofof{Proposition~\ref{prop: no fixpoint}}
Let $\calL$ be the set of partial injections without fixpoints (i.e.,
without size 1 cycles), let $L_n$ be the number of size $n$ elements
of $\calL$ and let $L(z)$ be the corresponding EGS. We want to show
that $L_n = \frac1e I_n (1+o(1))$.

The EGS $L(z)$ is computed using the standard
calculus of enumeration of labeled structures (displayed in
\cite[Figure II-18]{FlajoletSedgewick}), which was already used to
compute $I(z)$ in \cite{BNW}: since the EGS of cycles is
$\log(\frac1{1-z})$, the EGS of cycles of size at least 2 is
$\log(\frac1{1-z})-z$ and the EGS of non-empty sequences is
$\frac{z}{1-z}$, we have
$$L(z) = \exp\left(\log(\frac1{1-z}) - z + \frac{z}{1-z}\right) =
\frac1{1-z} \exp\left(\frac{z^2}{1-z}\right) = I(z) e^{-z}.$$
We already know that $I(z)$ is H-admissible \cite[Lemma 2.8]{BNW} and
Proposition~\ref{lemma: no fixpoint} shows that
$L(z)$ is H-admissible as well.

The saddlepoint is the solution $\zeta_n$ in the open interval $]0,1[)$
of the equation $z \frac{L'(z)}{L(z)} = n$.  An elementary computation
shows that we need to solve the equation
\begin{align*}
    &z^3 + (n-1)z^2 - (2n+1)z + n = 0,\\
    \textit{i.e.},\quad &(z+n+1)(1-z)^2 - 1 = 0.
\end{align*}
Letting $z = 0$ and $z = 1$ in this equation shows that there is a
solution in the interval $(0,1)$; moreover, one verifies easily that
$(z+n+1)(1-z)^2 - 1$ is monotonous on $(0,1)$, and hence our equation
has exactly one solution in that interval, say, $\zeta_n$.  From $0 < \zeta_n < 1$, we
deduce that $\frac1{n+2} < (1-\zeta_n)^2 < \frac1{n+1}$, and hence
$1-\sqrt{\frac1{n+1}} < \zeta_n < 1-\sqrt{\frac1{n+2}}$.  In particular,
$\zeta_n = 1 - \frac1{\sqrt n} + \O\left(\frac1{n\sqrt n}\right)$.

It now follows from Theorem~\ref{thcol} that
$$[z^n]L(z) = \frac{L(\zeta_n)}{\zeta_n^n\sqrt{2\pi b(\zeta_n)}}(1+o(1)).$$
In view of the proof of Proposition~\ref{lemma: no fixpoint}, $b(t) = b_0(t) - t$, where $b_0$ is the corresponding function for the H-admissible function $I(z)$. Using \cite[Equation~(7), p. 392]{BNW}, we find that $b_0(t) = \frac{2t}{(1-t)^3}$, and hence $b(t) = \frac{t(1+3t-3t^2+t^3)}{(1-t)^3}$.

Elementary computations show that
\begin{align*}
    b(\zeta_n) &= 2n^{\frac32}\left(1 + \O\left(\frac1{\sqrt n}\right)\right) \\
    \textrm{and }\frac1{\sqrt{2\pi b(\zeta_n)}} &= \frac{n^{-\frac34}}{2\sqrt{\pi}}
    \left(1 + \O\left(\frac1{\sqrt n}\right)\right). \\
    \textrm{Moreover, }\zeta_n^{-n} &= \exp\left(-n \log\left(1 - \frac1{\sqrt n} +
    \O\left(\frac1{n\sqrt n}\right)\right)\right) \\
    &= \exp\left(\sqrt n + \frac12 + \O\left(\frac1{\sqrt
    n}\right)\right). \\
    \textrm{Finally, }L(\zeta_n) &= \sqrt n \exp\left(\sqrt n - 2 + \O\left(\frac1{\sqrt
    n}\right)\right) \left(1 + \O\left(\frac1n\right)\right). \\
    \textrm{At last we have: }[z^n]L(z) &= \frac{e^{-\frac32}}{2\sqrt\pi}n^{-\frac14}e^{2\sqrt
    n}\left(1 + \O\left(\frac1{\sqrt n}\right)\right).
\end{align*}

Comparing with the estimate of $[z^n]I(z)$ in Proposition~\ref{prop:
facts graphs}, we find the announced result, namely
$$\frac{L_n}{I_n} = \frac{[z^n]L(z)}{[z^n]I(z)} = \frac1e(1+o(1)).$$
\eopo

\section{Finitely presented groups}\label{sec: fp groups}

One of the motivations for the study of subgroup distributions has
been the investigation of the statistical properties of
\textit{finitely presented groups}, see
\cite{1987Gromov,1991Champetier,1995Champetier,Olshanskii1992,1996ArjantsevaO,Arjantseva2000}.  Strictly
speaking, this would require a notion of distribution of these groups,
so that one would make a list of non-isomorphic groups and investigate
the frequency of groups with certain properties within that list.  No
such notion is available, as far as the authors are aware and current
literature operates rather with a notion of distribution of finite
{\em presentations}.

Recall that a \textit{finite presentation} is a pair $(A,R)$, where
$A$ is a finite set (the alphabet of generators) and $R$ is a tuple of
elements of $F(A)$ (the relators).  The resulting finitely presented
group $G$, written $G = \langle A \mid R \rangle$, is the quotient $G
= F(A) / N(R)$, where $N(R)$ is the \textit{normal subgroup generated}
by $R$.  The usual approach of statistical properties of finitely
presented groups is based on the uniform distribution on $k$-tuples of
reduced (or cyclically reduced) words of length at most $n$.

Of course, different presentations may yield the same group, even if
the alphabet of generators is fixed.  We are not aware of an analogue
of Proposition~\ref{distinct subgroups} above, which would state, say,
that the distribution of finitely presented groups with $k$ relators
of length at most $n$ resulting from the uniform distribution on
$k$-tuples of reduced (or cyclically reduced) words of length at most
$n$, is uniform, at least on a generic subset of $k$-tuples.  However,
partial results exist in this direction for one-relator groups
(Kapovich, Schupp and Shpilrain \cite{2006KSS}, Sapir and Spakulova 
\cite{2010Sapir,2009SapirSpakulova}).

In this section, we want to discuss an idea that may seem reasonable
in this context, but which turns out to be disappointing.  If $H$ is
the subgroup generated by the tuple of relators $R$, then $N(R) =
N(H)$, so the group $G = \langle A \mid R \rangle$ is also specified
by the pair $\langle A \mid H \rangle$.  Thus, instead of looking at
the normal closure of a finite set of elements, we look at the normal
closure of a finitely generated subgroup.  Now, clearly, if one
generates a list of subgroups $H$ by listing $k$-tuples of generators
(the word-based distribution discussed earlier in this article), then
the distribution of groups produced by this process will be the same
as if one were working with presentations.

The idea we wish to explore is to generate the subgroup $H$ via its
Stallings graph, that is, to use the graph-based distribution of subgroups.
Precisely, we may present groups via pairs, $\langle A \mid \Gamma
\rangle$ where $A$ is an alphabet and $\Gamma$ is a Stallings graph.
This is {\em a priori} a more compact representation of the group
(more compact in bit size, less convenient to \LaTeX).%
\footnote{%
It would be more interesting to have a unique, discrete
representation of finitely generated normal subgroups, but no such
representation seems to be known. And distinct normal subgroups may
well lead to isomorphic quotients.}
More importantly, as we have seen that the graph-based distribution of
subgroups is different from the word-based distribution, we may
anticipate a different distribution of finitely presented groups as
well, which would give us different insights on finitely presented
groups.

Now an interesting feature of the statistical study of group
presentations by tuples of relators is that the groups produced are
generically non-trivial, and in fact infinite.  More strongly, if $A$
and $k$ are fixed and if the maximal length $n$ of the relators in the
$k$-tuple $R$ tends to infinity, then generically $G = \langle A \mid
R\rangle$ is such that every subgroup generated by $|A|-1$ elements is
free~\cite{1996ArjantsevaO}.  It is also known that $G$ is generically
hyperbolic (Ol'shanski\u\i\ \cite{Olshanskii1992} and Champetier \cite{1991Champetier,1995Champetier}, proving a statement of
Gromov \cite{1987Gromov}).

In sharp contrast, and somewhat disappointingly, generically, a 
finitely presented group of the form $\langle A \mid \Gamma\rangle$ 
is trivial. 

\begin{theorem}\label{trivial fp group}
    Generically, the finitely presented group $\langle A \mid
    \Gamma\rangle$ is trivial.  In other words, generically, the
    normal closure of a randomly chosen subgroup of $F_r$ of size $n$,
    is $F_r$ itself.
\end{theorem}

The rest of Section~\ref{sec: fp groups} is devoted to the proof of
Theorem~\ref{trivial fp group}.  We note that if the lengths of the
cycles of the partial injection induced by letter $a$ in $\Gamma(H)$
are relatively prime, then $a$ belongs to the normal subgroup $N(H)$,
and hence $a = 1$ in $G = \langle A \mid H\rangle$.  Thus it suffices
to prove the following proposition.

\begin{proposition}\label{relatively prime partial injections}
    Generically, the lengths of the cycles of a size $n$ partial
    injection are relatively prime.
\end{proposition}

\begin{remark}
    Our proof that $\langle A \mid H\rangle$ is generically trivial
    relies on a rather rough upper bound: we show that generically
    with probability $1 - \O(n^{-\frac16})$, each letter $a$ is a
    product of conjugates of powers of $a$ in $H$.  We do not know
    whether $\langle A \mid H\rangle$ is exponentially generically
    trivial.  See Remark~\ref{remark purity exp negligible} for a
    similar situation.
\end{remark}

\subsection{The permutation case}

We start with the case of permutations, which is interesting in and of
itself.  Observe that if the lengths of the orbits of a permutation
are not relatively prime, then these lengths have a common prime
divisor $p$, which is in particular a divisor of $n$.  Let
$\calP^{(p)}_n$ be the set of size $n$ permutations in which all the
orbits have size a multiple of $p$.

\begin{lemma}\label{lemma: Ppn}
    Let $n\geq 2$ and let $p$ be a prime divisor of $n$.Then
    $$|\calP^{(p)}_n|\leq 2n!\,n^{\frac1p-1}$$
\end{lemma}

\proof
We fix $p$, so $n$ is of the form $n = mp$ and we proceed by induction
on $m$.  If $m = 1$, that is, $p = n$, then $|\calP^{(p)}_n|$ is the
number of size $n$ cycles, namely $(n-1)!$.  We now assume that $m >
1$.

We enumerate the elements of $\calP^{(p)}_n$ in terms of the size 
$kp$ of the orbit of 1: to determine such a permutation, one needs to 
select the other $kp-1$ elements of that orbit, select a cycle on 
these $kp$ elements, and select a permutation on the remaining elements, 
that is, an element of $\calP^{(p)}_{n-kp}$. Thus, using the 
convention that $|\calP^{(p)}_{0}| = 1$, we have
\begin{align*}
    |\calP^{(p)}_n| &= \sum_{k=1}^{m}{{n-1}\choose
    {kp-1}}(kp-1)!|\calP^{(p)}_{n-kp}| \\
    &= \sum_{k=1}^{m}\frac{(n-1)!}{(n-kp)!}|\calP^{(p)}_{n-kp}| \\
    &= (n-1)!\sum_{j=0}^{m-1}\frac{|\calP^{(p)}_{jp}|}{(jp)!}.
    \end{align*}
Isolating the term $j = 0$ and using the induction hypothesis, it
follows that
$$|\calP^{(p)}_n|  \leq  
(n-1)!\left(1+2\sum_{j=1}^{m-1}(jp)^{\frac1p-1}\right) = 
(n-1)!\left(1 + 2p^{\frac1p-1}\sum_{j=1}^{m-1}j^{\frac1p-1}\right).$$
Since the map $x\mapsto x^{\frac1p-1}$ is non-increasing on positive 
reals, we have
$$ (j+1)^{\frac1p-1}\leq \int_{j}^{j+1}x^{\frac1p-1}dx = 
p((j+1)^{\frac1p} - j^{\frac1p}).$$
Therefore, isolating the term $j = 1$,
\begin{align*}
    \sum_{j=1}^{m-1}j^{\frac1p-1} &= 1 + 
    \sum_{j=2}^{m-1}j^{\frac1p-1} \\
    & = 1 + \sum_{j=1}^{m-2}(j+1)^{\frac1p-1} \\
    & \le 1 + p((m-1)^{\frac1p}-1) \\
    & \le 1 + p(m^{\frac1p}-1) \enspace=\enspace
    p^{1-\frac1p}n^{\frac1p} - p + 1.
\end{align*}
Now we have
\begin{align*}
    |\calP^{(p)}_n| & \leq (n-1)!\left(1+2p^{\frac1p-1}
    \left(p^{1-\frac1p}n^{\frac1p}-p + 1\right)\right)\\
    & \leq (n-1)!\left(2n^{\frac1p}+1-2p^{\frac1p} +
    2p^{\frac1p-1}\right).
\end{align*}
Since $p\geq 2$, it holds
$$p^{\frac1p}-p^{\frac1p-1} = p^{\frac1p}\left(1-\frac1p\right) \geq
p^{\frac1p}\frac12\geq \frac12 $$
and hence
$$|\calP^{(p)}_n| \leq (n-1)!2n^{\frac1p} = 2n!n^{\frac1p-1},$$
which concludes the proof.
\eop

\begin{proposition}\label{perm}
    The probability that the lengths of the orbits of a size $n$
    permutation are not relatively prime is at most equal to $
    \frac2{\sqrt n} + 2n^{-\frac23}\log_3n$.
\end{proposition}

\proof
Let $Q_n$ be the set of size $n$ permutations for which the lengths of
the orbits are not relatively prime, and let $q_n = \frac{|Q_n|}{n!}$.

As we already observed, a common divisor of the lengths of the orbits
of a size $n$ permutation is also a divisor of $n$.  Therefore, if $n$
is prime, then $Q_n$ is the set of size $n$ cycles, so $|Q_n| = 
(n-1)!$, $q_n = \frac1n$ and we have the desired result.

If $n$ is not prime, then every size $n$ permutation in $Q_n$ is in
$\calP^{(p)}_n$ for some prime divisor $p$ of $n$.  These
sets are not pairwise disjoint, but the sum of their cardinalities is an upper
bound for $|Q_n|$.  For these values of $p$, $|\calP_n^{(p)}| \le
2n!n^{\frac1p-1}$ by Lemma~\ref{lemma: Ppn}.  Separating the case
$p=2$ from the cases $p\ge 3$, we find that $q_n \le \frac2{\sqrt n} +
2Dn^{-\frac23}$, where $D$ is the number of distinct odd prime
divisors of $n$.  Since $n \ge 3^D$, we have $D \le \log_3n$ and hence
$$q_n \le \frac2{\sqrt n} + 2Dn^{-\frac23} \le \frac2{\sqrt n} + 
2n^{-\frac23}\log_3n,$$
which concludes the proof.
\eop

\subsection{Proof of Proposition~\ref{relatively prime partial injections}}

Isolating the cycles in a size $n$ partial injection, reveals a
permutation (on a subset $X$ of $[n]$) and a fragmented permutation
(i.e., a cycle-less partial injection) on the complement of $X$.

The EGS $J(z) = \sum_n\frac{J_n}{n!}z^n$ of fragmented permutations
was discussed in Section~\ref{sec: malnormal graph}, where we noted in
particular that $J(z) = \exp\left(\frac z{1-z}\right)$.  Let us add
the following observation.

\begin{lemma}\label{Jn/n! increasing}
    The sequence $(J_n/n!)_{n> 0}$ is increasing.
\end{lemma}

\proof
Let $M_n = \frac{J_n}{n!}$, so that $J(z) = \sum_{n\ge 0}M_nz^n$.  The
equalities $\frac{d}{dz}J(z) = \frac1{(1-z)^2}J(z)$, and hence
$(1-z)^2\frac{d}{dz}J(z) = J(z)$, yield the following recurrence
relation, for all $n\ge 2$:
$$(n+1)M_{n+1} = (2n+1)M_n - (n-1)M_{n-1}.$$
It follows that, for all $n\ge 2$,
$$(n+1)(M_{n+1} - M_n)= nM_n - (n-1)M_{n-1} = n(M_n - M_{n-1}) +
M_{n-1}.$$
The result follows by induction since $M_1 = 1$ and $M_2 = \frac32$
(see for instance \cite[Section II.4.2]{FlajoletSedgewick}).
\eop

Specifying a size $n$ partial injection whose permutation part (the 
union of the cycles) has size $k$, amounts to choosing $k$ elements, 
choosing a permutation on these $k$ elements, and choosing a 
fragmented permutation on the remaining $n-k$ elements: the number 
of such partial injections is
$${n \choose k}k!J_{n-k} = n! \frac{J_{n-k}}{(n-k)!};$$
and the number of those in which the sizes of the cycles have a
non-trivial gcd is at most equal to
$$2 n! \frac{J_{n-k}}{(n-k)!} \left(\frac1{\sqrt k} + \frac{\log_3 
k}{k^{\frac23}}\right) $$
by Proposition~\ref{perm}. Moreover, summing the numbers of partial 
injections with permutation part of size $k$, we get
$$I_n = \sum_{k=0}^n n! \frac{J_{n-k}}{(n-k)!}.$$

We use these observations to show the following facts, which together
suffice to establish Proposition~\ref{relatively prime partial
injections}.

\begin{fact}\label{first fact}
    The proportion of size $n$ partial injections whose permutation
    part has size less than $n^{\frac13}$ is 
    $\O(n^{-\frac16})$.
\end{fact}

\begin{fact}\label{second fact}
    The proportion of size $n$ partial injections whose permutation
    part has size greater than $n^{\frac13}$ and for which the sizes
    of the cycles has a non-trivial gcd, is $\O(n^{-\frac16})$.
\end{fact}

\proofof{Fact~\ref{first fact}}
The proportion of size $n$ partial injections whose permutation part
has size less than $n^{\frac13}$ is
\begin{align*}
    \frac1{I_n}\sum_{k=0}^{\lfloor n^{\frac13}\rfloor} n!
    \frac{J_{n-k}}{(n-k)!} &\le
    (n^{\frac13}+1)\frac{J_n}{I_n}\quad\textrm{by Lemma~\ref{Jn/n!
    increasing}}\\
    &\le \O(n^{-\frac16}).
\end{align*}
The last inequality holds since $\frac{J_n}{I_n} = \O(n^{-\frac12})$ 
(compare the asymptotic equivalents of $\frac{J_n}{n!}$ given in 
Section~\ref{proofofpropK} and of $\frac{I_n}{n!}$ in 
Proposition~\ref{prop: facts graphs}).
\eopo

\proofof{Fact~\ref{second fact}}
Here we use Proposition~\ref{perm} and the fact that, for large enough
integers, we have $\frac1{\sqrt k} + \frac{\log_3 k}{k^{\frac23}} \le
\frac2{\sqrt k}$.  The number of size $n$ partial injections whose
permutation part has size greater than $n^{\frac13}$ and for which the
sizes of the cycles has a non-trivial gcd, is bounded above by
\begin{align*}
    \sum_{k=\lceil n^{\frac13}\rceil}^{n} 2n!  \frac{J_{n-k}}{(n-k)!}
    \left(\frac1{\sqrt k} + \frac{\log_3 k}{k^{\frac23}}\right) &\le
    4 n^{-\frac16} \sum_{k=\lceil
    n^{\frac13}\rceil}^{n} n!\frac{J_{n-k}}{(n-k)!} \\
    &\le 4n^{-\frac16} \sum_{k=0}^{n} n!\frac{J_{n-k}}{(n-k)!} = 4
    n^{-\frac16}I_n.
\end{align*}
Thus the proportion of these partial injections is at most
$4n^{-\frac16}$.
\eopo

\subsection*{Acknowledgements}

We want to thank Nicolas Pouyanne for fruitful discussions on the gcd
of the lengths of cycles in a random permutation. We also thank the anonymous referees for their insightful remarks, pointing us to important references and spotting two computational mistakes.

\bibliographystyle{plain}
{\small


\newcommand{\Ju}{Ju}\newcommand{\Th}{Th}\newcommand{\Ch}{Ch}\newcommand{\Yu}{Y%
u}\def\cprime{$'$} \def\cprime{$'$}


}

\end{document}